\documentstyle[12pt]{article}
\title{Almost closed interscribed polygons}
\author{
Yury Kroll \& Boris Mirman\\
\footnotesize
Department of Mathematics and Computer Science\\
\footnotesize
Suffolk University, Beacon Hill, Boston, MA 02114,\\
\footnotesize
E-mail: mirman@mcs.suffolk.edu}
\date{}

\begin{document}

\maketitle

\underline{{\bf Abstract}}\\
\\
\normalsize
We demonstrate a new approach to the computation of ratios of elliptic integrals. It turns out that almost closed polygons interscribed between two conics retain some of the properties of such closed polygons. We apply these retained properties to compute ratios of an incomplete elliptic integral over the complete one. This computation is based on an iterative procedure to determine the sequence of vertices of a polygon interscribed between two conics. Surprisingly, some iteration numbers are the denominators of the convergent fractions for the ratio of some elliptic integrals. The algorithm ensures high precision, is numerically stable and as fast as the arithmetic-geometric mean method, though not faster. Nonetheless, there are reasons to consider the proposed algorithm as it is quite different from other integration methods and may be applicable to other problems.\\
AMS classification: 14N15.\\
Keywords: Poncelet's theorem; Continued fraction; Dynamical system.

\section{Introduction}
\normalsize
As is commonly known, if a definite integral cannot be expressed through elementary functions, then the main way to calculate it is through a summation following from its definition. In some cases it is possible to avoid summation. One such case is the integration of irrationals which lead to elliptic integrals: an elliptic integral can be then computed using the arithmetic-geometric mean method. This paper presents a different method of computing incomplete elliptic integrals which avoids summation, but requires the knowledge of the corresponding complete integral. So for one value of the complete integral, we may determine a set of incomplete integrals. Possible applications of this approach are the integration of other functions, problems related to Blaschke products, and dynamical systems.\\ 
\\
The geometric figures that arise in the following discussion are figures in the real plane ${\bf R}^2$ or, if it is specifically mentioned, the complex plane ${\bf C}^2$. We interpret the points of ${\bf R}^2$ as complex numbers. A segment with ends $z$ and $w$ is denoted $[z,w]$. We consider here the unit circle ${\cal C}$ of ${\bf R}^2$ and a domain inside ${\cal C}$. In order to distinguish the disc ${\cal C}$ from its boundary, the latter is notated $\partial\cal C$.\\ 
\\
The words ``matrices" and ``operators in a complex Hilbert space" are used interchangeably. The scalar product 
of vectors $x$ and $y$ in a Hilbert space is denoted by $<x,y>$. A convex domain in ${\bf R}^2$ sometimes may be interpreted as the numerical range $W(T)=\left\{<Tx,x>:\:||x||=1\right\}$ of a matrix $T$. Such an interpretation helps find clarifying propositions and their proofs. For the properties of the numerical ranges of matrices, the reader may consult, for example, [17].
Any numerical range of a finite matrix is a closed convex domain. Unitarily equivalent matrices have the same numerical range. The numerical range of a unitary matrix is the polygon inscribed in $\cal C$ that is the convex hull of the eigenvalues of this unitary matrix. The boundary of $W(T)$ is written $\partial W(T)$. If $T$ is a contraction (i.e., $||T||\le1$), then $W(T)$ is inside ${\cal C}$. The main object of this study is a polygon inscribed in ${\cal C}$ and circumscribed about a domain inside ${\cal C}$. For brevity, we say that such a polygon is {\it interscribed} between the domain and $\partial{\cal C}$. An interscribed polygon is not necessarily closed.\\
\\
\underline{{\bf Definition}}. We call ``almost closed interscribed polygons" a sequence of $N_j$-sided polygons ${\cal P}_{N_j}\subset{\cal P}_{N_{j+1}}$, $j=1,2,\ldots$, interscribed between the domain and $\partial{\cal C}$ where the distance between the initial and last vertices of each ${\cal P}_{N_j}$ monotonically approaches zero.\\
\underline{{\bf Remark 1.}} The monotonic approach is necessary here for the application of the theory of continued\\
\hspace*{0.9in} fractions.\\
\\
A closed convex polygon inscribed in ${\cal C}$ and containing $W(T)$ (not necessarily circumscribed about $W(T)$) sometimes may be considered as the numerical range of a unitary dilation of $T$; that is, a unitary operator $U$ acting in a Hilbert spaces $H_1$, with $T$ acting in a Hilbert space $H\subset H_1$, and $T=P_HUP_H$, where $P_H$ is the orthogonal projection from $H_1$ to $H$. If $dim(H_1\ominus H)=k$, then $U$ is called a unitary $k$-dilation of $T$. In a series of papers, Gau and Wu [8 - 13, 39] and Mirman {\it et al.} [30 - 35], have considered the Poncelet case of $W(T)$ and a closed interscribed polygon $W(U)$ interpreted as the numerical range of a unitary $1$-dilation of $T$. In Section 2, we explore the geometry of circles, ellipses and the geometry of elliptic curves. A detailed consideration of conics is helpful for possible generalizations to the case when the domain is not a conic. Such domains and the appearance of the so-called attractive polygons are introduced in Section 3. In terms of dynamical systems, we consider there transformations $\partial{\cal C}\rightarrow\partial{\cal C}$ and the possible locations of the iterated vertices.\\
\\
If $U$ is a unitary $1$-dilation of $T$, then $\partial W(U):={\cal P}$ is the closed polygon interscribed between $\partial W(T)$ and $\partial{\cal C}$. The location of the tangent points ${\cal P}\bigcap W(T)$ must satisfy the following condition: let $z_k$, $k=0,\ldots,N-1$ be the vertices of the polygon ${\cal P}$ and $\zeta_k\in\partial W(T)$ be the tangent point on the side $[z_{k-1},z_{k}]$ of ${\cal P}$, $z_{N}=z_0$. Then 
\begin{equation}
\Pi_{k=1}^{N}\:\:{\frac{z_{k}-\zeta_k}{\zeta_k-z_{k-1}}}=1.
\end{equation}  
If this condition is not satisfied, then there is no $N\times N$ unitary $1$-dilation of $T$ (see, for example, [31]). If nevertheless ${\cal P}$ is closed and interscribed between $\partial W(T)$ and $\partial{\cal C}$, then at least one vertex of ${\cal P}$ is a folded eigenvalue of $U$. See Example 7 in Section 3.2. We consider in Section 3 also polygons formed by diagonals and sides of $\partial W(U)$ interscribed between $\partial W(T)$ and $\partial{\cal C}$.

\section{The case of conics}

\subsection{Preliminaries}

First let us assume that the conic is a circle ${\cal K}$ centered at $c$ with radius $r$. Consider a polygon interscribed between $\partial{\cal K}$ and $\partial{\cal C}$ which never closes regardless of how many sides of this polygon are constructed. Based on this polygon, we may select a sequence of almost closed interscribed polygons. These polygons retain some of the properties of closed polygons interscribed between two circles. An iterative procedure to determine the sequence of vertices of such almost closed polygons will be studied. It happens that the numbers of the sides of almost closed polygons (i.e., the corresponding iteration numbers) are the denominators of the convergent fractions for the ratio of some Legendrian elliptic integrals. Based on this fact, we develop an algorithm to calculate the ratio of these integrals.\\
\\
In the complex plane ${\bf C}^2$, the common points of the circumferences $\partial{\cal K}$ and $\partial{\cal C}$ have the same real abscissa 
\begin{equation}
I:=\frac{1+c^2-r^2}{2c},
\end{equation}
and the same ordinates $\pm\sqrt{1-I^2}$, which may be real or pure imaginary. All those circumferences having the same common points as have $\partial{\cal K}$ and $\partial{\cal C}$ form a so-called {\it pencil}. If a chord $[z,w]$ of the circle ${\cal C}$ is tangent to ${\cal K}$, then $z$ and $w$ satisfy the equation
\begin{equation}
p(z,w):=(c w z-w-z+c)^2-4r^2w z=0.
\end{equation}
See [31] for the derivation of this equation. Equation (3) is valid for any mutual locations of ${\cal C}$ and ${\cal K}$ (nested, intersected, or separated) where the tangent point may be between the chord's ends or in the continuation of the chord. Though the generalization does not require much effort, for convenience, we consider in detail only the case of nested circles ${\cal K}\subset{\cal C}$, $0<c<c+r<1$. Then, if $|z|=1$ and $p(z,w)=0$, we have $|w|=1$ also. Consequently, the tangent point is always between $z$ and $w$, and $I>1$. Since $p(z,w)=p(w,z)$, Equation (3) is easy to use for determining the sequence of vertices of a polygon interscribed between $\partial{\cal K}$ and $\partial{\cal C}$.\\
\\
\underline{{\bf The procedure to determine the sequence of vertices $z_k:=e^{i\varphi_k}$ of a polygon interscribed between}}\\
\underline{{\bf $\partial{\cal K}$ and $\partial{\cal C}$ for ${\cal K}$ to be to the left of the chord $[z_k,z_{k+1}]$:}}\\
Let $z_{-1}=e^{-i\varphi_1}$ and $z_0=1$, where $\cos\varphi_1=2r^2/(1-c)^2-1$, $\sin\varphi_1>0$. For $k=0,1,\ldots$, we have from Equation (3)
\begin{equation}
z_{k+1}=\frac{(c-z_k)^2}{z_{k-1}(1-cz_k)^2}.\Box
\end{equation}
\underline{{\bf Remark 2}.} Gau and Wu [8] gave the equation with a Blaschke product on the left side
\[z\Pi_{j=1}^{n}\frac{z-w_j}{1-\bar{w}_jz}=b\]
to determine the vertices $z=z_k$ for a general $n\times n$ UB-matrix (unitary bordering matrix, see the definition in [31]) with eigenvalues $w_j$ and any chosen $b\in\partial C$. Here, the procedure is to determine the solution of such an equation for the particular case of a conic, allowing $n=\infty$. Livshitz [29] first initiated the study of operators that are related to the UB-matrices.\\
\\
It may happen that for some $N$, $z_N=z_0$. This would mean that there exists an $N$-sided closed polygon interscribed between $\partial{\cal K}$ and $\partial{\cal C}$. This is the Poncelet case: regardless of the starting point $z_0$, the procedure results in $z_N=z_0$. Here we mainly consider interscribed polygons which never close.\\ 
\\
In order to develop the algorithms mentioned in the Abstract, we need to determine some special indices of the vertices produced by this procedure. We show how to do this in a very quick way below. For that, we combine some basic elements of the theory of continued fractions with the theory of elliptic curves. Equation (3) defines an elliptic curve, denoted here by ${\cal E}_{z,w}$. Indeed, substitute  
\[w=\frac{y+zr^2+c(1-2Iz+z^2)}{(1-c z)^2}\]
into Equation (3). Then we have the equation of the transformed elliptic curve ${\cal E}_{z,y}$, almost in Weierstrass form,
\begin{equation}
y^2=4c r^2(z^3-2Iz^2+z), \hspace*{0.25in} {\sf where} \hspace*{0.25in} y=w(1-c z)^2-z r^2-c(1-2I z+z^2).
\end{equation}
Since we apply both Equations (3) and (5) for derivations and calculations, it is convenient to have there the same variable $z$. As usual, the curve ${\cal E}_{z,y}$ is equipped with a group structure, i. e., with the addition rule (addition is denoted by $\circ$). The identity point $O_y=(\ldots,\infty)$ is the common point of all vertical lines, $z=const$. The curve ${\cal E}_{z,w}$ is equipped with a group structure also, assuming that the addition on ${\cal E}_{z,w}$ (denoted here by $\oplus$) is induced by the addition on ${\cal E}_{z,y}$. In other words, $(z_1,w_1)\oplus(z_2,w_2)=(z_3,w_3)$ iff $(z_1,y_1)\circ(z_2,y_2)=(z_3,y_3)$, where $y_1$, $y_2$ and $y_3$ are defined by the second equation of (5). 
The addition of points of ${\cal E}_{z,y}$ can be carried out as described in, for example, Silverman and Tate [38]. Namely, for the addition $(z,y)\circ(Z,Y)=(\tilde{z},\tilde{y})$, we have
\begin{equation}
\tilde{z}=\frac{\left(Slope\right)^2}{4c r^2}+2I-z-Z, \hspace*{0.25in} \tilde{y}=\left(Slope\right)(Z-\tilde{z})-Y,
\end{equation}
where $Slope=(Y-y)/(Z-z)$, if $z\neq Z$, and $Slope=4cr^2(3z^2-4Iz+1)/(2y)$, if $z=Z$ and $y=Y\neq0$.

\subsection{The algorithm}

Legendre's elliptic integral of the first kind is
\begin{equation}
F(\psi,k):=\int_{0}^{\psi}{\frac{d t}{\sqrt{1-k^2\sin^2t}}}.\\
\end{equation}
\\
In this section, we develop the algorithm to calculate the quotient of $\beta(\psi,k)=F(\psi,k)/F(\pi/2,k)$ (an incomplete integral $F(\psi,k)$, i. e., $0<\psi<\pi/2$ by the complete integral $F(\pi/2,k)$) for the so-called normal case ($0<k^2<1$). Such a ratio is used in problems with a probabilistic measure defined on $\partial{\cal C}$ (Kolodziej [27], King [24]). Another application of this algorithm is to compute the incomplete integrals $F(\psi,k)$.\\
\\
To facilitate the derivations, we consider an integral that differs slightly from $F(\psi,k)$. Substitute $t=\pi/2-\phi/2$ into Equation (7). Then 
\begin{equation}
F(\psi,k)=\frac{1}{k\sqrt{2}}\left[\Phi(\pi,I)-\Phi(\pi-2\psi,I),\right] \hspace*{0.25in}{\sf where}\hspace*{0.25in}\Phi(\varphi,I):=\int_{0}^{\varphi}{\frac{d\phi}{\sqrt{I-\cos\phi}}},\:\:I=\frac{2}{k^2}-1,
\end{equation}
and $\beta(\psi,k)=1-2\theta(\pi-2\psi,I)$, where $\theta(\varphi,I)=\Phi(\varphi,I)/\left(2\Phi(\pi,I)\right)$. We suppose that the given are the center $c$ and radius $r$ of the circle ${\cal K}$. If $\psi$ and $k^2$ are given, then a proper choice of $c$ and $r$ will make the expression for $I$ in Equation (7) coinciding with Equation (2). Namely,
\begin{equation} 
c=\frac{\left(\sqrt{1-k^2\cos^2\psi}-\sqrt{1-k^2}\right)^2}{k^2\sin^2\psi}\hspace*{0.1in}{\sf and}\hspace*{0.1in}r=(1-c)\cos\psi.\\
\end{equation}
\\
The value of the upper limit $\psi$ as a function of the value of the integral $F(\psi,k)$ is a Jacobian elliptic function. Jacobi, [21] 1828, observed a similarity between the trigonometric and elliptic functions. In terms of $\Phi(\varphi,I)$, a part of Jacobi's findings is the following proposition:\\
\\
\underline{{\bf Proposition 1.}} Consider a chord $[e^{i\varphi},e^{i\psi}]$ of ${\cal C}$ which is tangent to the circle ${\cal K}$. If the tangent point $\zeta=\zeta(\varphi,\psi)$ is between the chord's ends $e^{i\varphi}$ and $e^{i\psi}$, then the integral $\int_{\varphi}^{\psi}{\frac{d\phi}{\sqrt{I-\cos\phi}}}$ does not depend on the choice of $\varphi$, as long as the chord is tangent to ${\cal K}$. $\Box$\\
\\
For more details and for the case when ${\cal K}$ is a conic other than a circle, see Bos et al. [2], Schoenberg [37], Kolodziej [27], and King [24,25].
Proposition 1 allows us to interpret the determination of the vertices of an interscribed polygon as a walk on the unit circumference with a constant ``step length" defined by a measure. We show below how to find the so-called best approximations of this step length through the measure density function $h(x)$. In order to work with these terms, we need some basic notations, definitions and a proposition from the theory of continued fractions. Let us employ the following taken from Khinchin [23] and Cassels [3].\\
\\
\underline{{\bf Notations:}} for any real number $x$,\\  
$[x]$ is the greatest integer that is smaller than $x$;\\ 
$\{x\}=x-[x]$ is the fractional part of $x$;\\
$||x||=\min\left(\{x\},1-\{x\}\right)$ is the distance from $x$ to the nearest integer.\\
\underline{{\bf Definition.}} A fraction $p/q$ ($q>0$) is called {\it a best approximation} of $\alpha$, if
$||q\alpha||=|q\alpha-p|<||j\alpha||$ for all positive integers $j$ that are smaller than $q$.\\
\underline{{\bf Proposition 2.}}\\
A. For any real number $\alpha$, there exist sequences of integers $q_0=1<q_1<q_2<\ldots$ and $p_0,p_1,p_2,\ldots$ such that
the fractions $p_n/q_n$ (the {\it convergents}) are the best approximations of $\alpha$ (all of them).\\  
B. If $\alpha$ is rational, then $\alpha=p_N/q_N$ for some $N$.\\
C. If $\alpha$ is irrational, then $p_n/q_n\rightarrow\alpha$ for $n\rightarrow\infty$.\\
D. The approximated $\alpha$ is always between successive convergents: $(q_n\alpha-p_n)(q_{n+1}\alpha-p_{n+1})\le0$.\\
E. $q_n||q_{n+1}\alpha||+q_{n+1}||q_n\alpha||=1$.\\
F. The ratio $a_n=\left(q_{n+1}-q_{n-1}\right)/q_n$ is an integer that depends on $n$ and $\alpha$, and the numerators satisfy the equation $p_{n+1}=a_np_n+p_{n-1}$.\\
\\
\underline{{\bf Lemma.}} Let $h(x)$ be a periodic non-negative function with period 1, satisfying $h(1-\{x\})=h(\{x\})$ and $\int_0^1h(t)dt=1$. Consider a walk along the $x$-axis in the positive direction starting at the origin with a step of constant irrational length 
\[\alpha=\int_{x_{k-1}}^{x_k}h(t)dt,\] 
i. e., the first step is from 0 to $x_1$; the $k$th step is from $x_{k-1}$ to $x_k$; the total distance walked after the $n$th step is 
\[n\alpha=\int_{0}^{x_n}h(t)dt.\]   
If the numbers $q_j$, $j=1,2,\ldots$ are in ascending order and such that $||x_{q_j}||<||x_{k}||$ for all $1\le k\leq q_j-1$, then $q_j$ are the denominators of the best approximations of $\alpha$.\\ 
{\bf Proof}. Put $x_n=n\alpha$, $||x_{q_{j}}||=\varepsilon_j$, and $\int_0^{\varepsilon_j}h(t)dt=\delta_j$. Note that since $h(1-t)=h(t)$ for $0<t<1$, we have $\int_{\{x\}}^1h(t)dt=\int_0^{1-\{x\}}h(t)dt$. Since $n\alpha=\int_0^{x_n}h(t)dt$ and $\int_0^1h(t)dt=1$, we have $[n\alpha]=[x_n]$, $\{n\alpha\}=\int_0^{\{x_n\}}h(t)dt$. Thus 
\[||x_n||=\min\left(\int_0^{\{x_n\}}h(t)dt,\int_0^{1-\{x_n\}}h(t)dt\right).\\\] 
\\
By the choice of the numbers $q_j$, we have for all $n<q_j$, $||n\alpha||>||q_j\alpha||$, i. e., the $q_j$ are the denominators of the best approximations of $\alpha$. $\Box$\\
\\
\underline{{\bf Remark 3.}} The point of this lemma is that one can find the best approximations of $\alpha$ based on information about the upper limits $x_k$, regardless of the particular values of the function $h(t)$. Note that if there is no symmetry of $h(x)$ in $(0,1)$ with respect to $0.5$, then the lemma is not true. Indeed, consider the example of a piecewise constant function $h(x)=0$ if $0\le x<0.5$, and $h(x)=2$ if $0.5\le x<1$. For $\alpha<1/4$, the first four $q_j$ are $1,\:2,\:3,\:4$, and assertion F of Proposition 2 is not satisfied. Nevertheless, a more extensive analysis shows that without this symmetry, only some initial values of $q_j$ are unsuitable for being denominators of best approximations of $\alpha$. Then, after these values, there exist convergents $p_j/q_j$ which $\rightarrow\alpha$.\\ 
\\ 
Applying the lemma to the case of 
\[h(x)=\frac{\pi}{\Phi(\pi,I)\sqrt{I-\cos(2\pi x)}},\]
we have the following corollary.\\
\underline{{\bf Corollary.}} Consider the sequence of vertices $z_0=1$, $z_k=e^{i\varphi_k}$, $k=1,2,\ldots$ of a polygon interscribed between $\partial{\cal K}$ and $\partial{\cal C}$ by the procedure described in Section 2.1. Let $x_k=\varphi_k/(2\pi)$ and
$\int_{x_k}^{x_{k+1}}h(x)dx=\theta.$ If the numbers $q_j$, $j=1,2,\ldots$ are in ascending order and such that $||\cos\varphi_{q_j}||<||\cos\varphi_{k}||$ for all $1\le k\leq q_j-1$, then the $q_j$ are the denominators of the best approximations of $\theta$. $\Box$\\
\\
If $\theta$ is rational, the continued fraction is finite, and we have the Poncelet case of a closed polygon for any starting point $z_0\in\partial{\cal C}$. Cie$\check{s}$lak, Martini and Mozgawa [5] considered a ratio, similar to $\theta$, to determine the rotation index of these polygons. If $\theta$ is irrational, then its continued fraction is infinite and the convergents $p_j/q_j\rightarrow\theta$. The $z_k$ are the vertices of the almost closed $q_j$-sided interscribed polygons. The vertices form a dense set in $\partial{\cal C}$. Sufficient conditions for such a case are given in [32]. In particular, we have such a case if $c$ and $r$ are rational and the square roots $\sqrt{(1\pm r)^2-c^2}$ are irrational. Some extensions of these facts are treated in Section 3 (the regular case).\\
\\
\underline{{\bf Example 1:}} Let $c=0.5$, $r=0.2$, $I=1.21$. Then the procedure shows that the sequence  
\[{\bf q}=\left(q_j\right)_{j=1}^{14}=\left(2,5,7,12,31,43,74,117,191,308,1115,9228,56483,291643\right)\]
is such as the corollary claims, i. e., we have 
\[{\bf |z_{q}-1|}=\left(|e^{i\varphi_{{q_j}}}-1|\right)_{j=1}^{14}=\left(0.570,0.312,0.222,0.0837,0.0519,0.0317,0.0202,0.0115,0.00869,\right.\]
\[\left.0.00278,0.000341,0.0000552,0.00000978,0.00000112\right),\\\]
\\
and $z_{q_j}$ is closer to $z_0=1$ than any $z_k$ for $0<k<q_j$: $|z_k-1|>|z_{q_j}-1|$.\\
\\
It is possible to determine the numbers $q_j$ in a much faster way than using the procedure of Section 2.1. We call below the chords of ${\cal C}$ which are tangents to ${\cal K}$ just {\it tangential} chords, omitting the word ``circle" when this will not cause confusion.\\   
\\
Following Griffiths [15] and Griffiths and Harris [16], we may interpret Proposition 1 by applying the addition rule on ${\cal E}_{z,w}$. Namely, for any tangential neighboring chords $[z_0,z_1]$ and $[z_1,z_2]$ of ${\cal C}$ (as in the corollary), the difference of the points of ${\cal E}_{z,w}$, 
$(z_1,z_2)\ominus(z_0,z_1)=\left(Z,W\right)$, does not depend on the choice of $z_0$. Hence, $(z_1,z_2)=(z_0,z_1)\oplus\left(Z,W\right)$,
$(z_2,z_3)=(z_0,z_1)\oplus[2]\left(Z,W\right)$, and so on, 
\begin{equation}
(z_k,z_{k+1})=(z_0,z_1)\oplus [k]\left(Z,W\right), 
\end{equation}
where 
$[k](Z,W)=\left(Z_k,W_k\right)=\left(Z,W\right)\oplus\left(Z,W\right)\oplus\ldots\oplus\left(Z,W\right)$ ($k$ items).\\
\\
Consider the tangential chord $[z_0=1,z_1=e^{i\varphi_1}]$ and its neighboring tangential chord $[e^{i\varphi_1},z_2=e^{i\varphi_2}]$ ($\varphi_2\neq0$). The addition
$(z_0,z_1)\oplus(Z,W)=(z_1,z_2)$ of the points of ${\cal E}_{z,w}$, or the addition $(z_0,y_0)\circ(Z,Y)=(z_1,y_1)$ of the points of ${\cal E}_{z,y}$
and Equations (5) and (6) yields $y_0=i(1-c)^2\sin\varphi_1$ and 
\begin{equation}
Z=\frac{1}{c}, \hspace*{0.2in} W=\frac{(1-c^2)^2}{4c r^2}, \hspace*{0.2in} Y=-\frac{2r^2}{c}\\
\end{equation}
\\
Note that the addition $(z_0,y_0)\circ(Z_k,Y_k)=(e^{i\varphi_k},y_k)$ yields
\begin{equation}
\cos\varphi_k=1-\frac{4Z_k(I-1)}{(Z_k-1)^2}.
\end{equation}
If the point $\left(Z,W\right)$ is of finite order $N$, i. e., $[N](Z,W)=O_w$ (the identity point of ${\cal E}_{z,w}$), then the vertex $z_N$ coincides with the starting vertex $z_0$ and the $N$-sided polygon with vertices $z_k$ closes. If the point $\left(Z,W\right)$ is of infinite order, then, by the corollary, the choice of the numbers $q_j$ ($j=1,2,\ldots$) is such that $\cos\varphi_{q_j}\rightarrow1$ monotonically (i. e., if $q_j<q_k$, then $\cos\varphi_{q_j}<\cos\varphi_{q_k}$). It follows from Equation (12) that either $|Z_{q_j}|\rightarrow\infty$ or $Z_{q_j}\rightarrow0$. Below (inequalities (13) and Equation (15)) we will show that in the case where ${\cal K}\subset{\cal C}$, only $|Z_{q_j}|\rightarrow\infty$ is possible, and this approach to the limit is also monotonic. By the definition of the identity point $O_y$, we have $\left(Z_{q_j},Y_{q_j}\right)\rightarrow O_y$ and $\left(Z_{q_j},W_{q_j}\right)\rightarrow O_w$. This fact suggests considering other points of ${\cal E}_{z,w}$ for the determination of the indices $q_j$.\\
\\
\underline{{\bf The procedure to determine $w_k$ for $w_0=0$ such that $p(w_k,w_{k+1})=0$ for $k=1,2,\ldots$:}}\\ 
The $w_k$ constitute the sequence of solutions of Equation (3) starting with $z=w_0=0$:\\
$w_1=c$,\\
$w_2=4cr^2/(1-c^2)^2$, and for $k=2,3,\ldots$,
\begin{equation}
w_{k+1}=\frac{(c-w_k)^2}{w_{k-1}(1-cw_k)^2}.
 \Box \end{equation}

This procedure is correct since for any $k>0$, the denominator is not equal to zero. Indeed, if for some $k>0$, $w_k=0$, then the $k$-sided polygon closes (see for example [31]). This contradicts the assumed condition. Also it is easy to prove by induction that
\begin{equation}
0<w_k<I-\sqrt{I^2-1}\hspace*{0.2in}(<1<1/c).
\end{equation}
and
\begin{equation}
Z_k=\frac{1}{w_k}.\\
\end{equation}
\\
The numbers $w_k$ have a simple geometric meaning: they are the centers of the circles belonging to the pencil defined by ${\cal C}$ and ${\cal K}$ and having rational ratios of the integrals $\Phi(2\psi_k,I)$, where $\psi_k=\arcsin\left(\sqrt{1+w_k^2-2Iw_k}/(1-w_k)\right)$ (see [31] and Berger [1]).\\ 
\\
\underline{{\bf Remark 4.}} The numbers $w_k$ are the zeroes of the Blaschke product mentioned in Remark 2 above. In a series of papers, Chalendar,
Gorkin, Partington [4], as well as Daepp, Gorkin, Voss [6] solved related problems, in particular problems which link the zeroes of some Blaschke products with their values on $\partial C$.\\ \\
Equations (11) and (14) yield
\begin{equation}
\cos\varphi_k=1-\frac{4w_k(I-1)}{(w_k-1)^2}.\\
\end{equation}
\\
Put $\gamma_k=\sqrt{w_k/c}.$ Obviously, $\gamma_{q_j}\rightarrow0$ also monotonically. Thus, the numbers $q_j$ can be determined directly from Equation (3) using the following theorem.\\
\\
\underline{{\bf Theorem 1.}} Consider the following procedure for iteratively finding a sequence of solutions to Equation (3): $\gamma_1=1$,\\
 $\gamma_2=2r/(1-c^2)$, $\gamma_{k+1}=\left(1-\gamma_k^2c^2\right)/\left(\left|\gamma_k^2-1\right|\gamma_{k-1}\right)$, $k=2,3,\ldots$. 
Let $q_j$, $j=1,2,\ldots$ be the all positive integers $q_1<q_2<\ldots$ such that the corresponding positive numbers $\gamma_{q_1}>\gamma_{q_2}>\ldots\rightarrow 0$ and satisfy the inequalities $\gamma_k>\gamma_{q_j}$ for all $k<q_j$. Then $q_j$ are the denominators of the best approximations of $\theta=\Phi(\varphi,I)/\left(2\Phi(\pi,I)\right)$,
where $\Phi(\varphi,I)$ is the elliptic integral defined in Equations (8). $\Box$\\
\\
Continuing Example 1, we have\\
\[{\bf q}=\left(q_j\right)_{j=1}^{25}=\left(2,5,7,12,31,43,74,117,191,308,1115,9228,56483,291643,348126,1336021,1684147,\right.\]  
\[\left.6388462,14461071,237765598,252226669,489992267,1232211203,21437582718,2702367633671\right);\]
Based on the values of $q_j$, we can determine the numerators $p_j$ of the best approximations of $\theta$, applying assertion F of Proposition 2: 
$p_{n+1}=a_np_n+p_{n-1}$, where $a_n=\left(q_{n+1}-q_{n-1}\right)/q_n$.\\
\\
For Example 1, we have
\[{\bf p}=\left(p_j\right)_{j=1}^{20}=\left(1,2,3,5,13,18,31,49,80,129,467,3865,23657,122150,145807,559571,705378,\right.\]
\[\left.2675705,6056788,99584313,105641101,205225414,516091929,8978788207,1131843406011\right).\]   

The results of the approximation of $\theta$ are
\small
\[\frac{p_{12}}{q_{12}}=\frac{3865}{9228}=0.4188339835\]
\[\ldots\ldots\ldots\ldots\ldots\ldots\ldots\ldots\ldots\]
\[\frac{p_{24}}{q_{24}}=\frac{8978788207}{21437582718}=0.418833985394304193770062   \\\]
\[\ldots<\frac{p_{22}}{q_{22}}<\frac{p_{24}}{q_{24}}<\ldots<\theta=0.418833985394304193770083<\ldots<\frac{p_{25}}{q_{25}}<\frac{p_{23}}{q_{23}}<\ldots;\\\]
\[\frac{p_{25}}{q_{25}}=\frac{898394912629}{2144990483003}=0.418833985394304193770084   ;\]
\[\ldots\ldots\ldots\ldots\ldots\ldots\ldots\ldots\ldots\]
\[\frac{p_{13}}{q_{13}}=\frac{23657}{56483}=0.4188339854\\\]
.\\

In accordance with Proposition 2, the difference $q_{j+1}-q_j$ increases with increasing $j$. The well known Baby-Step Giant-Step algorithm (Shanks) gives a hint for how to reduce the number of operations in searching for $q_{j+1}$. In accordance with Proposition 2, F, we can search for a $q_{j+1}$ of the form $q_{j+1}=q_jN+q_{j-1}$, where $N$ is an unknown integer. The Shanks algorithm uses a linear form $a+bN$ to determine the required $a$ and $b$ when $N$ is given. Our case is simpler: $a$ and $b$ are known and we should find $N$.\\ 
\\
We determine $w_{(K+1)q_j+q_{j-1}}$ by adding points on the elliptic curve, i. e., by a ``giant step" (Equations (5) and (11)):
\begin{equation}
\left(w_{(K+1)q_j+q_{j-1}},y_{(K+1)q_j+q_{j-1}}\right)=\left(w_{K q_j+q_{j-1}},y_{K q_j+q_{j-1}}\right)\circ\left(Z_{q_j},Y_{q_j}\right),\hspace*{0.2in} K=0, 1,\ldots,N.
\end{equation}
Let
\[\theta=\frac{1}{a_1+\frac{1}{a_2+\frac{1}{a_3+\ldots}}}\]
and
\[\frac{p_j}{q_j}=\frac{1}{a_1+\frac{1}{a_2+\frac{1}{a_3+\ldots+\frac{1}{a_j}}}}.\\\]

A giant step provides us with the number $a_{j+1}=N_{step}^{(j+1)}$ by the condition $w_{a_{j+1}q_j+q_{j-1}}<w_{q_j}<w_{(a_{j+1}-1)q_j+q_{j-1}}$. It may happen that we ought to try to avoid the subtraction of close large numbers (since the positive numbers $w_{q_{j+1}}$ and $w_{q_j}$ are small and the precision of the code may not be sufficient to avoid the loss of significant digits). Fortunately, it is possible to subtract the close large numbers algebraically and derive the equations without being affected to number rounding. See Appendix I for the corresponding Equations (A1 - A9). We can see there that $N$ giant steps require about $100N$ operations instead of about $10Nq_j$ operations of baby steps by the procedure for $w_k$. When $q_j>>10$, the giant steps are useful for reducing the number of operations.\\
\\
Note that we can derive a more precise equation for $\theta$ than just the equation for convergent fractions $\theta\approx p_j/q_j$. Such a precise equation is based on the following theorem.\\
\\
\underline{{\bf Theorem 2.}} For sufficiently small $\gamma_{q_{j-1}}=\Delta$, we have
\begin{equation}
N_{step}^{(j+1)}=\left[\frac{\gamma_{q_{j-1}}}{\gamma_{q_j}}\right],
\end{equation}
\begin{equation}
q_{j+1}=q_{j-1}+N_{step}^{(j+1)}q_j,
\end{equation}
\begin{equation}
\gamma_{q_{j+1}}=\gamma_{q_{j-1}}-N_{step}^{(j+1)}\gamma_{q_j}+O(\Delta^3),
\end{equation}
and
\begin{equation}
\theta\approx\frac{\gamma_{q_{j-1}}p_j+\gamma_{q_j}p_{j-1}}{\gamma_{q_{j-1}}q_j+\gamma_{q_j}q_{j-1}}.
\end{equation}
{\bf Proof.} Indeed, Equations (A1 - A4) of Appendix I show that $\epsilon_{cur}$, $v$ and $\alpha$ are of magnitude $O(\Delta^2)$. Due to Equation (A5), one giant step results in $\gamma_{new}=\gamma_{cur}-\gamma_{q_j}+O(\Delta^3)$, i. e., the values of the $\gamma_{new}$ form approximately a descending arithmetic sequence with a difference of $-\gamma_{q_j}$. This fact yields Equations (18 - 20). It follows from these equations that the continued fraction for the ratio $\gamma_{q_{j}}/\gamma_{q_{j-1}}$ is
\[\frac{\gamma_{q_j}}{\gamma_{q_{j-1}}}=\frac{1}{N_{step}^{(j+1)}+\frac{1}{N_{step}^{(j+2)}+\ldots}}.\]
The last continued fraction is the tail of the continued fraction for $\theta$ after determination of $q_j$, i. e., Equation (21):
\[\theta=\frac{p_{j-1}+p_j\times\gamma_{q_{j-1}}/\gamma_{q_j}}{q_{j-1}+q_j\times\gamma_{q_{j-1}}/\gamma_{q_j}}+O(\Delta^2).\:\:\Box\\\]
\\
For practical purposes, it is useful to know when the magnitude of $\Delta$ may be considered as ``sufficiently small". Based on our numerous calculations, we may say that if we need only three or four significant decimal digits for $\theta$ then $\Delta$ should be $<0.1$. Another extreme case, when we need 100 correct digits, requires that $\Delta$'s magnitude should be in the range of $10^{-24}$ - $10^{-26}$. Still the calculations are very fast and numerically stable: 100 significant decimal digits for $\theta$ may be obtained in less than 0.5 seconds by Phyton on today's consumer PC. The algorithm is described in detail in Appendix I.\\
\\
\underline{{\bf Remark 5.}} The giant step procedure may be interpreted in the following way: instead of iterating from side to side of the polygon, we iterate from diagonal to diagonal of this polygon, a diagonal that leaps over a certain number of vertices.

\subsection{The case of other elliptic integrals}

Though other elliptic integrals may be transformed so as to fall under the cases of Sections 2.1 and 2.2 by a substitution, it may be convenient to calculate these integrals directly by an algorithm which is similar to that considered above. Let   
\[\theta=\int_0^{\psi_1}\frac{d\varphi}{\sqrt{\alpha_0-2\alpha_1\cos\varphi+\alpha_2\cos^2\varphi}}/\int_0^{2\pi}\frac{d\varphi}{\sqrt{\alpha_0-2\alpha_1\cos\varphi+\alpha_2\cos^2\varphi}}.\]  
Then, for suitable $\alpha_0$, $\alpha_1$, $\alpha_2$ and $\psi_1$, we can find an ellipse $(\xi-c)^2/a^2+\eta^2/b^2=1$ with $a\ge b>0$ and $c$ such that $a^2(1-b^2)+b^2c^2=\alpha_0$, $b^2c=\alpha_1$, $b^2-a^2=\alpha_2$ and
\[\frac{a^2+b^2-(1-c)^2}{(1-c)^2+b^2-a^2}=\cos\psi_1.\]
Hence, instead of a circle, we have an ellipse, denoted here by ${\cal E}$, and instead of Equation (3), we have (see [32]) the following equation for the ends of a tangential chord $[z,w]$:
\[w^2\left[(cz-1)^2+(b^2-a^2)z^2\right]-2w\left[(cz-1)(z-c)+(a^2+b^2)z\right]+(z-c)^2+b^2-a^2=0. \hspace*{1.5in} \mbox{(3-ell)}\]
Here, as also in Sections 2.1 and 2.2, the calculation of the sequence of vertices of a polygon interscribed between $\partial{\cal E}$ and $\partial{\cal C}$ allows us to determine the ratio $\theta$.\\
\\
Let $z_k=e^{i\psi_k}$ be a sequence of solutions of Equation (3-ell) starting with $z_0=1$ and with $z_1=\cos\psi_1+i\sin\psi_1$ ($\sin\psi_1>0$). For $k=1,2,\ldots$,
\[z_{k+1}=\frac{1}{z_{k-1}}\:\frac{(c-z_k)^2+b^2-a^2}{(b^2-a^2)z_k^2+(cz_k-1)^2}.\]
It is well known (see for example King [24]) that either $z_k$, $k=1, 2, \dots$ are the vertices of a closed polygon interscribed between $\partial{\cal E}$ and $\partial{\cal C}$, or they form a dense set in $\partial{\cal C}$. Consequently, if a closed interscribed polygon does not exist, then there does exist a sequence of $q_j$-sided, $j=1, 2, \ldots$ almost closed polygons, $j=1, 2, \ldots$. The denominators $q_{j+1}$ of the convergents are determined as follows:\\ 
$q_0=1$, and $\varepsilon_0=1-\cos\psi_1$. Let $q_j$ and $\varepsilon_j$ be known. For the smallest $k>q_j$ such that $1-\cos\psi_k<\varepsilon_j$, we have $q_{j+1}=k$ and $\varepsilon_{j+1}=1-\cos\psi_k$. It is clear that all the ways mentioned in Section 2.2 to speed up these calculations are applicable here also.\\  
\\
\underline{{\bf Example 2.}} Let $a=0.5$, $b=0.4$, $c=0.4$, $\cos\psi_1=5/27$, $\sin\psi_1=8\sqrt{11}/27$; $\alpha_0=0.2356,$ $\alpha_1=0.064,$ $\alpha_2=-0.09,$ $\cos\psi_1=5/27$. Then the convergents $p_j/q_j$ are 
\begin{center}
\begin{tabular}{||c|c|c|c|c|c|c|c|c|c|c|c|c|c|c||}\hline\hline
$j=$ & -1 & 0 & 1 & 2 & 3 & 4 & 5 & 6 & 7 & 8 & 9 & 10 & 11 & 12\\ \hline
$q_j=$ & 0 &1 & 3 & 13 & 16 & 45 & 151 & 196 & 1327 & 12139 & 25605 & 37744 & 214325 & 252069\\ \hline
$a_j=(q_{j}-q_{j-2})/q_{j-1}$& & &3&  4  &   1  &   2   &   3   &   1   &   6    &    9   & 2 & 1 & 5 & 1\\ \hline
$p_j=a_jp_{j-1}+p_{j-2}$& 1 &0 & 1 & 4 &   5 & 14 &  47 &  61 &  413 &  3778 &  7969 & 11747 &  66704 &  78451\\ \hline
\hline
\end{tabular}
\end{center}
\[\frac{78451}{252069}<\theta=\frac{1}{2}\int_0^{\psi_1}\frac{d\varphi}{\sqrt{\alpha_0-2\alpha_1\cos\varphi+\alpha_2\cos^2\varphi}}/\int_0^{\pi}\frac{d\varphi}{\sqrt{\alpha_0-2\alpha_1\cos\varphi+\alpha_2\cos^2\varphi}}<\frac{66704}{214325}.\\\]
\\
\underline{{\bf Remark 6.}} An ellipse ${\cal E}$ with foci $f_1$, $f_2$ and minor axis $b$ may be interpreted (Donoghue [7]) as the numerical range of the $2\times2$ matrix \[T=\left(\begin{array}{cc}f_1&2b\\ 0&f_2\\ \end{array}\right).\]
It may happen that $||T||>1$ although $W(T)={\cal E}$ is inside ${\cal C}$ (as it is in Example 2). This may seem to obstruct the application of the unitary dilation of $T$. However it is possible to find a finite matrix $\hat{T}$ such that $W(\hat{T})={\cal E}$ and has norm not bigger than $1$. More about matrices with coinciding numerical ranges can be found in Helton and Spitkovsky [20].\\ 
\underline{{\bf Remark 7.}} Consider a set of confocal ellipses with the foci $c+ie$ and $c-ie$ being inside ${\cal C}$. If the minor axis $a\rightarrow0$, then the corresponding ellipses approach to a ``degenerate ellipse," the segment ${\cal K}_0=\left[c-ie,c+ie\right]$. For this segment, there exists an ``interscribed degenerate polygon," the chord ${\cal P}_0=\left[c-i\sqrt{1-c^2},c+i\sqrt{1-c^2}\right]$ of ${\cal C}$, and we may create a set of iterated vertices $z_k=e^{i\psi_k}$ as above. However, unlike the case of a real ellipse, this set of vertices is neither finite nor dense in $\partial{\cal C}$. The set of vertices is broken in two subsets where the vertices have two accumulation points, the ends of the chord ${\cal P}_0$. In the next section, we will see more examples of $N$-sided closed interscribed polygons whose vertices $Z_j$, $j=1,\dots,N$, are the accumulation points for the vertices $z_k=e^{i\psi_k}$.

\section{Possible generalizations}

Here we consider the case when the domain ${\cal K}$ is not a conic. There are the definitions, basic assertions and various examples. This section may be considered as an introduction to the further analysis of the interscribed polygons. We assume in this section that ${\cal K}$ is symmetrical with respect to the horizontal axis. The domain is strongly inside ${\cal C}$. The border $\partial{\cal K}$ is a smooth curve and may contain a flat portion. The border is assumed to be generated by a square matrix $T$, namely by an eigenvalue $\lambda(\phi)$ of the matrix $\Re\left(e^{-i\phi}T\right)=\left(e^{-i\phi}T+e^{i\phi}T^{*}\right)/2$, where the entries of $T$ are real. The tangent points $\zeta$ of $\partial{\cal K}$ are
\begin{equation}
\zeta=\xi+i\eta=\left(\lambda+i\lambda^{'}_{\phi}\right)e^{i\phi},
\end{equation}
i. e.,
\[\xi=\lambda\cos\phi-\lambda^{'}_{\phi}\sin\phi, \hspace*{0.25in} \eta=\lambda\sin\phi+\lambda^{'}_{\phi}\cos\phi.\]
If $\partial{\cal K}$ contains a flat portion, then some additional analysis is required to determine the points $\zeta$ (see Equation (29) with $a=b$ below and Rodman and Spitkovsky [36]). The curve $\partial{\cal K}$ has a simple relation with the $k$-numerical range $W_k(T)$ (Halmos [18,19], Li et al. [28], Gau et al. [14]). If $\lambda(\phi)$ is the largest eigenvalue of $\Re\left(Te^{-i\phi}\right)$, then ${\cal K}$ is the numerical range $W(T)$. Toeplitz noticed this in 1918 introducing the numerical range. Then Hausdorff proved that there are no holes inside $W(T)$. In 1951, Kippenhahn [26] proposed an equation for $\partial W(T)$ which is similar to Equation (22). We consider here the case when $\lambda(\phi)$ is not necessarily the largest eigenvalue of $\Re\left(Te^{-i\phi}\right)$, and we still apply Equation (22) to the points of $\partial{\cal K}$. Then the domain ${\cal K}$ may be not convex. However, the border $\partial{\cal K}$ satisfies the following (see [31]):\\
$\bullet$ For any angle $\phi$, $0\le\phi<2\pi$, there exists exactly one directed tangent line to $\partial{\cal K}$ that forms the angle $\phi$ with
\hspace*{0.1in} the vertical axis.\\
This property is satisfied for any convex domain, but $\partial{\cal K}$ may contain cusps (see for example [34]). This property ensures that for any point $\tilde{z}$ outside ${\cal K}$, there exists a line containing $\tilde{z}$ and tangent to $\partial{\cal K}$. This allows us to develop a simple iteration procedure to calculate the sequence of vertices of polygons interscribed between $\partial{\cal K}$ and $\partial{\cal C}$.\\
\\
By this procedure, we can also obtain, as for conics, the following interscribed polygons:
\begin{enumerate}
\item closed polygons with one arbitrarily chosen vertex on $\partial{\cal C}$ (Poncelet's case); 
\item almost closed polygons with a dense set of vertices in $\partial{\cal C}$ for one starting vertex.\\ 
\\ 
However, unlike conics, there is one other case possible: 
\item there exist only one or two closed interscribed polygons.
\end{enumerate}
Though these three cases do not exhaust all possibilities for interscribed polygons, below we consider only the last two cases. Let the chord $[z=e^{i\varphi},w=e^{i\psi}]$ of ${\cal C}$ be tangent to $\partial{\cal K}$ at the point $\zeta$, and if we move along the chord from $z$ to $w$, then ${\cal K}$ is to our left. Then (see [35])
\begin{equation}
\det\left(T+wzT^{*}-(w+z)I\right)=0.
\end{equation}
\begin{equation}
\frac{d\psi}{d\varphi}=\frac{e^{i\psi}-\zeta}{\zeta-e^{i\varphi}}
\end{equation}
Since $\zeta$ is between $e^{i\psi}$ and $e^{i\varphi}$, the derivative in Equation (24) is positive. Moreover, this derivative lies between two positive numbers, because the domain ${\cal K}$ is strongly inside ${\cal C}$ and the absolute values of the numerator and denominator on the right side of (24) are separated from zero by a positive number. For a given $z$, we can find $w$ from Equation (23). In other words, Equation (23) determines the sequence of vertices $z_k=e^{i\psi_k}$ of a polygon interscribed between $\partial{\cal K}$ and $\partial{\cal C}$. If $\partial{\cal K}$ does not contain any flat portion, then the point of tangency $\zeta_k$ on the chord $[z_{k-1},z_k]$ is unique, and we can define a function $h\left(z\right)$ on the set $\{z_k=e^{i\psi_k}\}_{k=0}^{\infty}\subset\partial{\cal C}$ as follows. Let us choose an arbitrary positive value $h(z_0)$ for the starting point $z_0=e^{i\psi_0}$, and for $k=1, 2, \ldots$, put

\[h(z_{k})=h({z_{k-1}})\:\frac{d\psi_{k}}{d\psi_{k-1}}=h({z_{k-1}})\:\frac{e^{i\psi_k}-\zeta_k}{\zeta_k-e^{i\psi_{k-1}}}\]

and
\begin{equation}
h(z_{k})=h({z_{0}})\:\frac{d\psi_{k}}{d\psi_{0}}=h(z_0)\Pi_{k=1}^{n}{\frac{z_{k}-\zeta_k}{\zeta_k-z_{k-1}}}.
\end{equation}
The behavior of the function $h(z_k)$ depends on the last product, $P_n=\Pi_{k=1}^{n}{\frac{z_{k}-\zeta_k}{\zeta_k-z_{k-1}}}.$ (See Figures 1a and 1b.)

\subsection{Regular interscribed polygons}

As stated in the Introduction, if $T$ is a UB-matrix (the Poncelet case) of size $N-1$, then for any starting point $z_0$, we have a closed $N$-sided polygon interscribed between the $\partial W(T)$ and $\partial{\cal C}$, and in accordance with Equation (1), the product $P_{N}=1$. Naturally, for almost closed interscribed $q_j$-sided  polygons, we have $P_{q_j}\rightarrow1$. More precisely, we have the following proposition.\\ 
\\
\underline{{\bf Proposition 3}}. Suppose that $\partial{\cal K}$ does not contain any flat portion. If the sequence $z_k$, $k=0,1,\ldots$ is dense in $\partial{\cal C}$, and the values of the ratios $\left|h(z_k)-h(z_l)\right|/|z_k-z_l|$ for any $k\neq l$ are bounded between two positive numbers, then $h(z_k)$ may be continued to a differentiable function $h(z)$ defined on all of $\partial{\cal C}$ with a bounded derivative. There exists an infinite sequence of numbers $q_j$ such that $z_k$, $k=1,\ldots,q_j$ are the vertices of the $q_j$-sided almost closed polygons, and the products
\[P_{q_j}=\Pi_{k=1}^{q_j}{\frac{z_{k}-\zeta_k}{\zeta_k-z_{k-1}}}\rightarrow1, \hspace*{0.2in} \mbox{when} \hspace*{0.2in} j\rightarrow\infty.\] 
{\bf Proof.} Let $\hat{z}\in\partial{\cal C}$. Because of the given conditions, there exists a subsequence $z_{k_j}$ such that $z_{k_j}\rightarrow\hat{z}$ and the $\lim h(z_{k_j})$ exists. This limit is $h(\hat{z})$. Without loss of generality, we may assume that $|z_{k_j}-\hat{z}|\rightarrow0$ monotonically. If $\hat{z}$ is the starting point $z_0$, then the $k_j=q_j$ are the numbers of sides of the almost closed polygons. Since $h\left(z_{q_j}\right)\rightarrow h(z_0)$, Equation (25) yields $P_{q_j}\rightarrow1$. $\Box$\\ 
\\
Let us call a polygon {\it regular} if its vertices $z_k$ satisfy Proposition 3. The features of such vertices and the corresponding continued functions $h(z)$ are formulated in the following assertions:
\begin{enumerate}
\item  There is no closed polygon interscribed between $\partial{\cal K}$ and $\partial{\cal C}$.    
\item For different starting points $z_0^{(1)}$ and $z_0^{(2)}$, the continued functions $h^{(1)}(z)$ and $h^{(2)}(z)$ are constant multiples of each other.
\item For any tangential chords $[s_1=e^{i\varphi_1},t_1=e^{i\psi_1}]$ and $[s_2=e^{i\varphi_2},t_2=e^{i\psi_2}]$, the integral $J(s_1,s_2)=\int_{\varphi_1}^{\varphi_2}d\varphi/h(s)$ along the $arc(s_1,s_2)$ of $\partial{\cal C}$ is equal to $J(t_1,t_2)=\int_{\psi_1}^{\psi_2}d\psi/h(t)$. In particular, the integral $J_0=J(s,t)$ along the $arc(s,t)$ spanned by a tangential chord $[s=e^{i\varphi},t=e^{i\psi}]$ is a constant, i. e., the measure of the $arc(s,t)$ is constant for the measure density $1/h$.
\item The numbers $q_j$ of Proposition 3 satisfy the condition for being the denominators of convergents, i. e., the ratio $(q_{j+1}-q_{j-1})/q_j$ is an integer (at least for some $j>j_0$, see Remark 3 after the proof of the lemma in Section 2.2). The corresponding convergents $p_j/q_j\rightarrow J_0/J(0,2\pi)$.
\end{enumerate}
{\bf Proof.} Assertions 1 and 2 are obvious. Assertion 3 follows from the definition of $h(z)$, i. e., Equation (24). Indeed, for any tangential chord $[s=e^{i\varphi},t=e^{i\psi}]$, we have $d\psi/h\left(t\right)=d\varphi/h\left(s\right)$; therefore for any pair of tangential chords $[s_1,t_1]$ and $[s_2,t_2]$, we have 
\[\int_{\varphi_1}^{\varphi_2}\frac{d\varphi}{h\left(s\right)}=\int_{\psi_1}^{\psi_2}\frac{d\psi}{h\left(t\right)}.\]
Assertion 4 follows from Assertion 3 and the lemma of Section 2. $\Box$\\
\\
In order to demonstrate the features of regular polygons, we consider below $3\times3$ matrices of upper triangular form 
\begin{equation}
T=\left(\begin{array}{ccc}c_1 & b_1 & a\\ 0 & c_2 & b_2\\ 0 & 0 & c_3\\ \end{array}\right).\\
\end{equation}
Applying Equation (23) to the matrix of (26) with $b_1=b_2=b$ and $c_1=c_2=c_3=0$, we have  
\begin{equation}
w^3+w^2\left[(3-A)z-B z^2\right]+w\left[(3-A)z^2-B z\right]+z^3=0,
\end{equation}
where $A=a^2+2b^2$ and $B=ab^2$. A detailed description of $W(T)$ for a $3\times3$ matrix $T$ is given by Keeler, Rodman and Spitkovsky [22] and by Rodman and Spitkovsky [36]. In particular, it was established there that for $a\neq b,\:\:ab\neq0$, $\partial W(T)$ is an oval not containing a flat portion. It is known also that for $a=b\neq0$, $W(T)$ is the convex hull of a cardioid, i. e., it is a cardioid together with a flat portion on the boundary, whereas for $ab=0$, $W(T)$ is a circle. Below, Examples 3 - 6 are for regular polygons.\\  
\\
\underline{{\bf Example 3.}} Let $a=0.6$, $b=0.4$. The convergent fractions $p_j/q_j$ are:
\begin{center}
\small
\begin{tabular}{||c|c|c|c|c|c|c|c|c|c|c|c|c|c|c|c|}\hline\hline
$j=$ &  -1 & 0 & 1 & 2 & 3 & 4  & 5  &   6 &   7 &   8 &     9 & 10 & 11 & 12 & 13\\ \hline
$q_j=$ & 0 &1  & 2 & 3 & 8 & 11 & 19 & 182 & 201 & 383 & 10925 & 11308 & 78773&247627&1564535\\ \hline
$(q_{j}-q_{j-2})/q_{j-1}$& & &2&  1  &   2  &   1 &   1   &   9   &   1   &    1   & 28 & 1&6&3&6\\ \hline
$p_j=a_jp_{j-1}+p_{j-2}$& 1 &0 & 1 & 1 &3& 4 & 7&  67 &  74 & 141&4022&4163&29000&91163&575978\\ \hline\hline
\end{tabular}
\end{center}
\begin{center}
\small
\begin{tabular}{|c|c|c|c|c|c|c|c|c||}\hline\hline
$j=$   & 14     & 15 & 16 & 17 & 18 & 19 & 20 & 21 \\ \hline
$q_j=$ &6505767&27587603&89268576&116856179&206124755&322980934&529105689&852086623\\ \hline
$(q_{j}-q_{j-2})/q_{j-1}$ &   4   &      4 &      3 &       1 &       1 &       1 &       1 &        1\\ \hline
$p_j=$ &2395075&10156278&32863909&43020187 &75884096 &118904283&194788379&313692662\\ \hline\hline
\end{tabular}
\end{center}
\[\frac{194788379}{529105689}<\frac{J_0}{J(0,2\pi)}<\frac{313692662}{852086623}.\\\]

\underline{{\bf Example 4.}} Let $a=0.72000001$, $b=0.72$. The convergent fractions $p_j/q_j$ are:
\begin{center}
\small
\begin{tabular}{||c|c|c|c|c|c|c|c|c|c|c|c|c|}\hline\hline
$j=$ &  -1 & 0 & 1 & 2 & 3 & 4  & 5  &   6 &   7 &   8 &     9 & 10\\ \hline
$q_j=$ & 0 &1  & 3 & 7 & 24 & 103 & 1363 & 2829 & 9850 & 32379 & 106987 & 139366\\ \hline
$(q_{j}-q_{j-2})/q_{j-1}$& & &3&  2  &   3  &   4 &   13   &   2   &   3   &    3   & 3 & 1\\ \hline
$p_j=a_jp_{j-1}+p_{j-2}$& 1 &0 & 1 & 2 &7& 30 & 397& 824 & 2869&9431&31162&40593\\ \hline\hline
\end{tabular}
\end{center}
\begin{center}
\small
\begin{tabular}{|c|c|c|c|c|c|c|c|c||}\hline\hline
$j=$ &  11 & 12 & 13 & 14 & 15 & 16 & 17 & 18\\ \hline
$q_j=$ &246343&632072&1510497&3653066&19775827&23428893&66633613&489864184\\ \hline
$(q_{j}-q_{j-2})/q_{j-1}$& 1 & 2 & 2&2&5&1&2&7\\ \hline
$p_j=$&71755&184103&439961&1064025&5760086&6824111&19408308&142682267\\ \hline\hline
\end{tabular}
\end{center}
\[\frac{142682267}{489864184}<\frac{J_0}{J(0,2\pi)}<\frac{19408308}{66633613}.\]
If $k\leq10^7$, we have $h_{max}/h_{min}\approx151$. If $k\leq10^8$, we have $h_{max}/h_{min}\approx227$. There is no high growth of $h_{max}/h_{min}$ when $k\rightarrow\infty$. Compare with Examples 9 and 10 in the next section.\\  
\\
\underline{{\bf Example 5.}} Let $a=0.2$, $b=0.4$, $c_1=c_3=0.1$, $c_2=0.35$ in Equation (27). The convergents $p_j/q_j$ are
\begin{center}
\small
\begin{tabular}{||c|c|c|c|c|c|c|c|c|c|c|c|c|}\hline\hline
$j=$ &  -1 & 0 & 1 & 2 & 3 & 4 & 5 & 6 & 7 & 8 & 9 & 10 \\ \hline
$q_j=$ & 0 & 1 & 2 & 3 & 5 & 58 & 179 & 416 & 2259 & 13970 & 44169 & 58139\\ \hline
$(q_{j}-q_{j-2})/q_{j-1}$&   &  &2&  1  &  1  & 11 & 3 &  2 & 5 & 6 &3 & 1\\ \hline
$p_j=a_jp_{j-1}+p_{j-2}$ & 1 &0 &1&  1  & 2 & 23 & 71 & 165 & 896 & 5541 & 17519 & 23060\\ \hline\hline
\end{tabular}
\end{center}
\begin{center}
\small
\begin{tabular}{|c|c|c|c|c||}\hline\hline
$j=$ &  11 & 12 & 13 & 14\\ \hline
$q_j=$ &102308 &160447 &262755&423202\\ \hline
$(q_{j}-q_{j-2})/q_{j-1}$& 1 & 1 & 1 & 1\\ \hline
$p_j=$ & 40579 & 63639 & 104218 & 167857\\ \hline\hline
\end{tabular}
\end{center}
\[\frac{167857}{423202}<\frac{J_0}{J(0,2\pi)}<\frac{104218}{262755}.\]

\normalsize
\underline{{\bf Example 6.}} Let $a=1-b^2=0.618033974844<b=0.618034$. The convergents $p_j/q_j$ are:

\begin{center}
\small
\begin{tabular}{||c|c|c|c|c|c|c|c|c|c|c|}\hline\hline
$j=$ &  -1 & 0 & 1 & 2 & 3 & 4  & 5  &   6 &   7 &   8 \\ \hline
$q_j=$ & 0 &1  & 3 & 274 & 6579 & 125275 & 257129 & 896662 & 11017073 & 22930808 \\ \hline
$(q_{j}-q_{j-2})/q_{j-1}$& & &3&  91  &   24  &   19 &   2   &   3   &   12   &    2  \\ \hline
$p_j=a_jp_{j-1}+p_{j-2}$& 1 &0 & 1 & 91 &2185& 41606 & 85397&  297797 & 3658961 & 7615719\\ \hline\hline
\end{tabular}
\end{center}
\begin{center}
\small
\begin{tabular}{|c|c|c|c|c|c||}\hline\hline
$j=$ & 9 & 10 & 11 & 12 & 13\\ \hline
$q_j=$ & 33947881 & 294513856 & 328461737 &622975593&951437330\\ \hline
$(q_{j}-q_{j-2})/q_{j-1}$& 1 & 8 & 1 & 1 & 1\\ \hline
$p_j=$& 11274680&97813159&109087839&206900998&315988837\\ \hline\hline
\end{tabular}
\end{center}
\[\frac{206900998}{622975593}<\frac{J_0}{J(0,2\pi)}<\frac{315988837}{951437330}.\]

\subsection{Attractive interscribed polygons}

Now let us consider the case when there exists a closed $N$-sided polygon ${\cal P}$ interscribed between $\partial{\cal K}$ and $\partial{\cal C}$ with vertices $Z_l=e^{i\varphi_l}$, $l=0,1,2,\ldots,N$, $Z_N=Z_0$. If ${\cal P}$ is not symmetrical with respect to the horizontal axis, then the polygon $\bar{{\cal P}}$ with the vertices $\bar{Z}_l$ is also a closed interscribed polygon. {\it We assume in this section that} ${\cal P}$ and $\bar{{\cal P}}$ {\it are the only closed interscribed polygons.}\\ 
\\
Let the numbering of $Z_l$ be such that the chords $[Z_{l-1},Z_{l}]$, $Z_0=Z_N$, are tangent to $\partial{\cal K}$ and the domain ${\cal K}$ is to the left of such a chord. The tangent points of the chords $[Z_{l-1},Z_{l}]$  are denoted by $\zeta_l^{(\infty)}$, and
\begin{equation}
P_{N}=\Pi_{l=1}^{N}{\frac{Z_{l}-\zeta_l^{(\infty)}}{\zeta_l^{(\infty)}-Z_{l-1}}}.
\end{equation}
As was mentioned above, $P_N$ is positive. If $\delta$ is a small deviation from $\varphi_0$, then it follows from Equation (24) that
\begin{equation}
\varphi_N=\varphi_0+P_N\delta+\alpha\delta^2+O(\delta^3).
\end{equation}
Therefore
\begin{enumerate}
\item If $P_N<1$, then the deviation of $\varphi_N$ from $\varphi_0$ is smaller than $\delta$. 
\item If $P_N>1$, then the deviation of $\varphi_N$ from $\varphi_0$ is greater than $\delta$. 
\item If $P_N=1$, then the deviation of $\varphi_N$ from $\varphi_0$ is smaller than $\delta$ for $\alpha\delta<0$ and greater than $\delta$ for $\alpha\delta>0$.
\end{enumerate}
Let us define a transform ${\bf R}:\partial{\cal C}\rightarrow\partial{\cal C}$ by the condition ${\bf R}e^{i\varphi}=e^{i\psi}$, where 
$\psi>\varphi$ and chord $[e^{i\varphi},e^{i\psi}]$ is tangent to $\partial{\cal K}$. Then we may say that if $P_N<1$, then all vertices of ${\cal P}$ are attractors. If
$P_N>1$, then all vertices of ${\cal P}$ are repelling points. If $P_N=1$, then all vertices of ${\cal P}$ are from one side attractors, and from the other side are repelling points.\\ 
\\
Notice that
\[\Pi_{l=1}^{N}{\frac{\bar{Z}_{l}-\bar{\zeta}_{l-1}^{(\infty)}}{\bar{\zeta}_{l-1}^{(\infty)}-\bar{Z}_{l-1}}}=\Pi_{l=1}^{N}{\frac{Z_{l}-\zeta_{l-1}^{(\infty)}}{\zeta_{l-1}^{(\infty)}-Z_{l-1}}}=\frac{1}{P_N}.\]
Therefore, for a vertex $e^{-i\varphi_0}$ of $\bar{{\cal P}}$, we have instead of Equation (29),
\[\varphi_N=\varphi_0+\frac{1}{P_N}\delta-\frac{\alpha}{P_N^2}\delta^2+O(\delta^3).\]
Hence if the vertices of ${\cal P}$ are attractors then the vertices of $\bar{{\cal P}}$ are repelling points. Consequently, we call ${\cal P}$ the {\it attractive polygon} and $\bar{\cal P}$ the {\it repelling polygon}. Remind that there is no other closed interscribed polygons. Therefore, if $P_N\neq1$, then each attractor is between two neighboring repelling points on $\partial{\cal C}$. Then the points $z_k$ can be broken up into $N$ mutually disjoint subsets such that the $z_k$ of each subset tend to the corresponding attractor.\\
\\ 
Similarly to the regular case, we can define the function $h(z_k)$ by applying Equation (24). However, here the conditions of Proposition 3 are not satisfied: the sequence $z_k$ is not dense in $\partial{\cal C}$ and the ratios $\left|h(z_k)-h(z_l)\right|/|z_k-z_l|$ can be unbounded or  arbitrarily close to zero. Consequently, Assertions 1-4 are not valid for this case, and the product $P_{N}$ of Equation (28) can be any positive number. Hence $h(z_k)$ may be arbitrarily big or arbitrarily close to zero.
Moreover, we have the following sufficient condition for an attractive/repelling polygon:\\
\underline{{\bf Proposition 4}}. If $P_N\neq1$, then $h(z_k)$ either exponentially grows to $\infty$ or exponentially falls to zero.\\
This proposition may be applied in order to distinguish the regular case from the attractive/repelling one, when other methods to distinguish the cases are more difficult to apply. Below, Examples 7-10 are for attractive polygons.\\
\\
\normalsize
\underline{{\bf Example 7.}} Let $a=1-b^2=0.618035210911>b=0.618033$. Consider the unitary two-dilation of $T$
\[U=\left(\begin{array}{ccccc}0 & b & a & 0 & -b\sqrt{a}\\ 0 & 0 & b & 0 & \sqrt{a}\\ 0 & 0 & 0 & 1 & 0\\ 1 & 0 & 0 & 0 & 0\\
 0 & \sqrt{a} & -b\sqrt{a} & 0 & b^2\\ \end{array}\right).\]
The eigenvalues of $U$ are $-1,1,1,-a/2\pm i\sqrt{1-a^2/4}$. The triangle $\Delta$ with vertices $(1,0)$, $(-a/2,\pm\sqrt{1-a^2/4})$ is the attractive triangle for Example 6 and is not such a triangle for Example 5. Indeed, for $a=1-b^2<b$, we have regular cases, because $\Delta$ is not interscribed between $\partial W(T)$ and $\partial{\cal C}$: its vertical side is not tangent to the curve $\partial W(T)$. For $a=1-b^2>b$, in contrast, the vertical side of $\Delta$ is tangent to $\partial W(T)$, and $\Delta$ is the attractive triangle (see Figure 2). The threshold for these cases is $1-b^2=b=(\sqrt{5}-1)/2$, when we have the cardioid $\partial W(T)$, which combines the curves $\partial{\cal K}_1$and $\partial{\cal K}_2$. Notice that the triangle $\Delta$ is symmetrical with respect to the horizontal axis. Therefore, if $a=1-b^2>b$, we have the case $P_N=1$, and the vertices of $\Delta$ behave as attractors from one side and repellers from another side.\\
\\
\underline{{\bf Example 8.}} Let $a=0.21$, $b_1=b_2=0.20$, $c_2=0.66$ and $c_1=c_3=0$ in Equation (27).\\ 
Then there is the attractive pentagon (see Figure 3) with the vertices\\
$(-0.997219,0.074522)$, $(0.938000,-0.346636)$, $(0.045972,0.998943)$, $(-0.253912,0.967227)$, $(0.970625,0.240598)$.\\
\\
\underline{{\bf Example 9.}} Let $a=0.7200001$, $b=0.72$. Then we have the 18337-sided attractive polygon. $P_{18337}=0.7029723633$. The starting point
$z_0=0.997910504956172999592891236-i*0.064611331035011368320516583$ yields a cycle of length 18337, with the convergents $p_j/q_j$ presented in the table
\begin{center}
\small
\begin{tabular}{||c|c|c|c|c|c|c|c|c|c|c||}\hline\hline
$j=$ &  -1 & 0 & 1 & 2  & 3 & 4  & 5  &   6&   7&   8 \\ \hline
$q_j=$ & 0 &1  & 3 & 7  &24 & 103&1363&2829&15508&18337\\ \hline
$(q_{j}-q_{j-2})/q_{j-1}$& & &3&  2&3& 4 & 13&  2 & 5  & 1\\ \hline
$p_j=a_jp_{j-1}+p_{j-2}$     &1&0&1& 2 &7& 30&397&824&4517&5341\\ \hline\hline
\end{tabular}
\end{center}
For $k\leq10^7$, $h_{max}/h_{min}\approx10^{85}$. For $k\leq10^8$, $h_{max}/h_{min}\approx10^{837}$.\\ 
\\
\underline{{\bf Example 10.}} Let $a=b=0.72$. Then we have the 3750742-sided attractive polygon, $P_{3750742}=0.6852390384$. The starting point
$z_0=0.715565891923305685013680-i*0.698545241423$ yields a cycle of length 3750742, with the convergents $p_j/q_j$ presented in the table
\begin{center}
\small
\begin{tabular}{||c|c|c|c|c|c|c|c|c|c|c|c|c|}\hline\hline
$j=$ &  -1 & 0 & 1 & 2 & 3 & 4  & 5  &   6 &   7 &   8 &     9 & 10\\ \hline
$q_j=$ & 0 &1  & 3 & 7 & 24 & 103 & 1363 & 2829 & 9850 & 71779 & 81629 & 153408\\ \hline
$(q_{j}-q_{j-2})/q_{j-1}$& & &3&  2  &   3  &   4 &   13   &   2   &   3   &    3   & 3 & 1\\ \hline
$p_j=a_jp_{j-1}+p_{j-2}$& 1 &0 & 1 & 2 &7& 30 & 397& 824 & 2869&9431&31162&40593\\ \hline\hline
\end{tabular}
\end{center}
\begin{center}
\small
\begin{tabular}{|c|c|c|c|c|c||}\hline\hline
$j=$ &  11 & 12 & 13 & 14 & 15 \\ \hline
$q_j=$ &235037&1093556&1328593&2422149&3750742\\ \hline
$(q_{j}-q_{j-2})/q_{j-1}$& 1 & 4 &1&1&1\\ \hline
$p_j=a_jp_{j-1}+p_{j-2}$&68459&318519&386978&705497&1092475\\ \hline\hline
\end{tabular}
\end{center}
For $k\leq10^7$, $h_{max}/h_{min}\approx195$. For $k\leq10^8$, $h_{max}/h_{min}\approx0.904\times10^{6}$. For $k\leq10^9$, $h_{max}/h_{min}\approx0.747\times10^{46}$.

\section{Appendix I}

\subsection{Determination of $\theta$ for circles (Section 2.2)}

\footnotesize
\begin{enumerate}
\item   
For given $\psi$ and $k^2$, $0<\psi<\pi/2$, $0<k^2<1$, calculate $c$ and $r$ using Equations (9).
\item
Set $j=1$, $q_1=1$, $\gamma_1=1$, and $\gamma_2=2r/(1-c^2)$; if $\gamma_2>1$, then $\varepsilon=1$ and $y_{q_1}=2cr^2$;\\
if $\gamma_2<1$, then $q_1=2$, $\varepsilon=\gamma_2$, and $y_{q_1}=-c\varepsilon^2((1-c^2)(1-\varepsilon^2)+2r^2)$.
\item
For $k=2, 3, \ldots$, calculate $\gamma_{k+1}=\left|1-\gamma_k^2\right|/((1-c^2\gamma_k^2)\gamma_{k-1})$. 
\item
If $k>2$ and $\gamma_k<\varepsilon$, then $j=j+1$, $q_j=k$, $\varepsilon=\gamma_k$ and\\
$y_{q_j}=c\gamma_{k+1}^2(1-c^2\varepsilon^2)-c\varepsilon^2r^2-c(1-2Ic\varepsilon^2+c^2\varepsilon^4)$.
\item
If $\varepsilon>0.1$, go to Item 3.
\item
Otherwise, go to the giant steps, Equation (16) and Equations (A1 - A9) of the next section of the Appendix I, with the parameters $c$, $r$, and $I$, and the two triples $q_{j-1}$, $\gamma_{q_{j-1}}$, $y_{q_{j-1}}$ and $q_{j}$, $\gamma_{q_{j}}$, $y_{q_{j}}$. Then calculate the triple $q_{j+1}$, $\gamma_{q_{j+1}}$, $y_{q_{j+1}}$.
\item
The final value of $\theta$ is calculated by Equation (21).
\end{enumerate}

\subsection{One set of the giant steps}

Input: constants $c$, $r$, $I=(1+c^2-r^2)/(2c)$, where $0<c<c+r<1$;\\
\\
\hspace*{0.4in} parameters which do not change during this one set of the giant steps:\\
\hspace*{1.4in} $q_j$, $\gamma_{q_j}$, $y_{q_j}$,\\
\hspace*{1.4in} $\varepsilon=\gamma_{q_j}$, $W_{q_j}=1/(c\varepsilon^2)$, $Y_{q_j}=-y_{q_j}W_{q_j}^2$,\\
\hspace*{1.4in} $\rho_{q_j}=\sqrt{1-2Ic\gamma_{q_j}^2+c^2\gamma_{q_j}^4}$,\\ 
\hspace*{1.4in} $\epsilon_{q_j}=1-\rho_{q_j}=c\gamma_{q_j}^2(2I-c\gamma_{q_j}^2)/(1+\rho_{q_j})$.\\                 
\\
\hspace*{0.4in} Initial variables for the first giant step:\\
\hspace*{0.4in} $q_{cur}=q_{j-1}$, $\gamma_{cur}=\gamma_{q_{j-1}}$, $y_{cur}=y_{q_{j-1}}$.\\
\\
(*) Then we calculate:

$\begin{array}{lc}
\rho_{cur}=\sqrt{1-2Ic\gamma_{cur}^2+c^2\gamma_{cur}^4}, & \hspace*{1.in}(A1)\\
\epsilon_{cur}=1-\rho_{cur}=c\gamma_{cur}^2(2I-c\gamma_{cur}^2)/(1+\rho_{cur}), & \hspace*{1.in} (A2)\\
v=4Ic\gamma_{cur}\gamma_{q_j}-2(\epsilon_{cur}+\epsilon_{q_j}-\epsilon_{cur}\epsilon_{q_j})-c^2\gamma_{cur}\gamma_{q_j}(\gamma_{cur}^2+\gamma_{q_j}^2), & \hspace*{1.in} (A3)\\
\alpha=\frac{\gamma_{cur}\gamma_{q_j}}{(\gamma_{cur}-\gamma_{q_j})^2}v, & \hspace*{1.in} (A4)\\
\gamma_{new}=\frac{\gamma_{cur}-\gamma_{q_j}}{1-c^2\gamma_{cur}^2\gamma_{q_j}^2}\sqrt{1-\alpha}, & \hspace*{1.in} (A5)\\
y_{new}=\frac{y_{cur}c\gamma_{q_j}^2+y_{q_j}}{1-c^2\gamma_{cur}^2\gamma_{q_j}^2}\:\frac{\gamma_{new}^2-\gamma_{cur}^2}{\gamma_{q_j}^2}-y_{cur}, & \hspace*{1.in} (A6)\\
\gamma_{cur}=\gamma_{new}, & \hspace*{1.in} (A7)\\
y_{cur}=y_{new}, & \hspace*{1.in} (A8)\\
q_{cur}=q_{cur}+q_j. & \hspace*{1.in} (A9)\\
\end{array}$\\

If $\gamma_{cur}>\varepsilon$, go to (*) and calculate (A1 - A9) again. If $\gamma_{cur}<\varepsilon$, then this set of the giant steps is completed.\\
For the next set of the giant steps, we have the following:\\
$q_{j+1}=q_{cur}$, $\gamma_{q_{j+1}}=\gamma_{new}$, $y_{q_{j+1}}=y_{new}$, $q_{cur}=q_j$, $\gamma_{cur}=\gamma_{q_j}$, $y_{cur}=y_{q_j}$.

\section{Appendix II}

\underline{{\bf A procedure to determine the sequence of vertices $z_k=e^{i\psi_k}$ of a polygon interscribed between $\partial W(T)$}}\\
\underline{{of Equation (27) and $\partial{\cal C}$ for $W(T)$ to be to the left of the chord $[z_k,z_{k+1}]$}:}\\
\\
\underline{Start}.\\
Given: $a$, $b_1$, $b_2$, $c_1$, $c_2$, $c_3$. Then $\varepsilon=10$; $j=-1$; $k=q_{-1}=0$;\\
$\alpha_1=c_1+c_2+c_3$,\\
$\alpha_2=c_1c_2+c_2c_3+c_3c_1$,\\
$\alpha_3=(a^2+b_1^2+b_2^2)/4$,\\
$\alpha_4=c_1c_2c_3$,\\
$\alpha_5=(c_1b_2^2+c_2a^2+c_3b_1^2-ab_1b_2)/4$.\\
\#3. $\cos\psi_0=1$, $\sin\psi_0=0$, $\lambda_1^2=(\alpha_3-\alpha_5)/(1-\alpha_1+\alpha_2-\alpha_4)$;\\
\#4. $\cos\psi_0=-1$, $\sin\psi_0=0$, $\lambda_1^2=(\alpha_3+\alpha_5)/(1+\alpha_1+\alpha_2+\alpha_4)$;\\
\#5. $\cos\psi_0=\sqrt{1-\alpha_3}$, $\sin\psi_0=\sqrt{\alpha_3}$, $\lambda_1^2=\alpha_3$.\\
\underline{*Iterations: k=1,2,...}.\\
$\lambda_{k}=\sqrt{\lambda_{k}^2}$;\\
$\cos\psi_k=(2\lambda_k^2-1)\cos\psi_{k-1}-2\lambda_k\sqrt{1-\lambda_k^2}\sin\psi_{k-1}$;\\
$\sin\psi_k=(2\lambda_k^2-1)\sin\psi_{k-1}+2\lambda_k\sqrt{1-\lambda_k^2}\cos\psi_{k-1}$;\\ 
$\beta_1=1-\alpha_1\cos\psi_k+\alpha_2(2\cos^2\psi_k-1)-\alpha_4(4\cos^3\psi_k-3\cos\psi_k)$;\\
$\beta_2=\alpha_5\cos\psi_k-\alpha_3+\alpha_2\sin^2\psi_k-3\alpha_4\cos\psi_k\sin^2\psi_k$;\\
$\beta_3=\left(\alpha_4-\alpha_1+2\alpha_2\cos\psi_k-4\alpha_4\cos^2\psi_k\right)\sin\psi_k$;\\
$\beta_4=(\alpha_5-\alpha_4\sin^2\psi_k)\sin\psi_k$.
\[p=\frac{\beta_3^2/2-\beta_1\beta_2-\beta_3\beta_4}{\beta_1^2+\beta_3^2}-\frac{\lambda_k^2}{2},\:\:q=\frac{\beta^2_4}{(\beta_1^2+\beta_3^2)\lambda_k^2}.\]
$\lambda_{k+1}^2=p+\sqrt{p^2-q}$.\\
\#3. $\delta=1-\cos(\psi_k)$.\\
\#4. $\delta=1+\cos(\psi_k)$.\\
\#5. $\delta=1-\cos(\psi_k-\psi_0)$.\\ 
If $\delta<\varepsilon$ then goto (*). If $\delta\geq\varepsilon$, then $\varepsilon=\delta$, $j:=j+1$, $q_j=k$, and then goto (*).\\

\Large
{\bf References}\\

\scriptsize

[1] M. Berger, {\it Geometry,} {\bf 2,} Springer Verlag, Berlin, 1987.\\

[2] H. J. M. Bos, C. Kers, F. Oort, D. W. Raven, {\it Poncelet's Closure Theorem,} Expo. Math. {\bf 5} (1987) 289-364.\\

[3] J.W.S. Cassels, {\it An Introduction to Diophantine approximation,} Cambridge University Press, (1961).\\

[4] I. Chalendar, P. Gorkin, J.R. Partington, {\it Determination of Inner Functions by their Value Sets on the Circle,} Computational Methods and Function Theory, 11 (2011), No. 1, 353-373.\\

[5] W. Cie$\check{s}$lak, H. Martini, W. Mozgawa, {\it On the rotation index of bar billiards and Poncelet's porism,} Bull. Belg. Math. Soc. {\bf 2} (2013) 287-300.\\  

[6] U. Daepp, P. Gorkin, K. Voss, {\it Poncelet's theorem, Sendov's conjecture, and Blaschke products,} J. Math. Anal. Appl. {\bf 365} (2010) 93-102.\\

[7] W. F. Donoghue Jr., {\it On the numerical range of a bounded operator,} Michigan Math. J. {\bf 4} (1957) 261-263.\\

[8] H.-L. Gau, P. Y. Wu, {\it Numerical range of $S(\phi)$,} Linear and Multilinear Algebra, 45 (1998) 49-73.\\

[9] H.-L. Gau, P. Y. Wu, {\it Lucas theorem refined,} Linear and Multilinear Algebra, 45 (1998) 359-373.\\

[10] H.-L. Gau, P. Y. Wu, {\it Numerical range and Poncelet property,} Taiwanese Journal of Math., 7(2) (Jun. 2003) 173-193.\\  

[11] H.-L. Gau, P. Y. Wu, {\it Condition for the numerical range to contain an elliptic disc,} Linear  Algebra Appl., 364 (2003) 213-222.\\ 

[12] H.-L. Gau, P. Y. Wu, {\it Numerical range of a normal compression,} Linear and Multilinear Algebra, 52 (2004) 195-201.\\

[13] H.-L. Gau, P. Y. Wu, {\it Numerical range circumscribed by two polygons,} Linear  Algebra Appl., 382 (2004) 155-170.\\

[14] H.-L. Gau, C. K. Li, P. Y. Wu, {\it Higher-rank numerical ranges and dilations,} Operator Theory 63 (2010) 181-189.\\  

[15] P. Griffiths, {\it Variations on a Theorem of Abel,} Inventiones Math. {\bf 35} (1976) 278-293.\\

[16] P. Griffiths and J. Harris, {\it A Poncelet theorem in space,} Math. Helv. {\bf 52} (1977) 145-166.\\

[17] K. E. Gustafson, D. K. M. Rao, {\it Numerical Range}, Springer Verlag, 2013.\\

[18] P. R. Halmos, {\it Numerical ranges and normal dilations,} Acta Szeged 25 (1964) 1-5.\\ 

[19] P. R. Halmos, {\it A Hilbert Space Problem Book,} P. 211, Springer-Verlag, Berlin, 1982.\\ 

[20] J. W. Helton and I. M. Spitkovsky, {\it The possible shapes of numerical ranges}, http://arxiv.org/pdf/1104.4587.pdf\\

[21] C.G.J. Jacobi, {\it $\ddot{U}$ber die Anwendung der elliptischen Transcendenten auf ein becanntes Problem der Elementargeometrie,} Journal $f\ddot{u}$r die reine and 
angewandte Mathematik ({\it Crell's Journal}) {\bf 3} (1828) 376-389.\\

[22] D. S. Keeler, L. Rodman, I. M. Spitkovsky, {\it The numerical range of $3 \times 3$ matrices,} Lin. Algebra Appl. {bf 252} (1997) 115-139.\\

[23] A. Ya. Khinchin, {\it Continued Fractions,} Dover Publications (1964).\\

[24] J. L. King, {\it Three problems in search of a measure,} Amer. Math. Monthly {\bf 101} (1994) 609-628.\\

[25] J. L. King, {\it Billiards Inside a Cusp,} Math. Intelligencer {\bf 17} (1995) 8-16.\\

[26] R. Kippenhahn, {\it $\ddot{U}$ber den Wertevorrat einer Matrix,} Math. Nachr. {\bf 6} (1951) 193-228.\\

[27] R. Kolodziej, {\it The rotation number of some transformation related to a billiards in an ellipse,} Studia Matematica {\bf LXXXL} (1985) 293-302.\\

[28] C. K. Li, C. H. Sung, and N. K. Tsing, {\it c-Convex matrices: Characterizations, inclusion relations and normality,} Lin, Multilin. Alg. 25 (1989) 275-287.

[29] M. S. Livshitz, {\it On a certain class of linear operators in Hilbert space,} Mat. Sbor. {\bf 19} (1946) 209-262 (in Russian with English resume).\\ 

[30] B. Mirman, {\it Numerical ranges and Poncelet curves,} Lin. Algebra Appl. {\bf 281} (1998) 59-85.\\

[31] B. Mirman, {\it UB-matrices and conditions for Poncelet polygon to be closed,} Linear Algebra Appl. {\bf 360} (2003) 123-150.\\

[32] B. Mirman, {\it Sufficient conditions for Poncelet polygons not to close,} Amer. Math. Monthly {\bf 112} (2005) 351-356.\\

[33] B. Mirman, {\it Explicit solutions to Poncelet's porism,} Linear Algebra Appl. {\bf 436} (2012) 3531-3552.\\

[34] B.Mirman, V.Borovikov, L.Ladyzhensky, R.Vinograd, {\it Numerical ranges, Poncelet curves, and invariant measures,} Linear Algebra Appl. {\bf 329} (2001) 61-75.\\

[35] B. Mirman, P. Shukla, {\it A characterization of complex plane Poncelet curves,} Linear Algebra Appl. 408 (2005) 86-119.\\

[36] L. Rodman, I. M. Spitkovsky, {\it $3 \times 3$ matrices with a flat portion on the boundary of the numerical range,} {\bf 397} (2005) 193-207.\\

[37] I. J. Schoenberg, {\it On Jacobi-Bertrand's proof of a Theorem of Poncelet,} Studies in Pure Mathematics, Birkh$\ddot{\mbox{a}}$user, Boston (1987) 623-627.\\

[38] J. H. Silverman, J. Tate, {\it Rational Points on Elliptic Curves,} Springer-Verlag, New York, 1992.\\

[39] P. Y. Wu, {\it Polygons and numerical ranges,} Amer. Math. Monthly {\bf 107} (2000) 528-540.

\end{document}